\documentclass[12pt]{article}
\usepackage{amsmath}
\usepackage{amssymb}
\usepackage{amsfonts}
\usepackage{eucal}
\usepackage{graphicx}
\numberwithin{equation}{section}
\newtheorem{Theorem}{Theorem}[section]
\newtheorem{Proposition}[Theorem]{Proposition}
\newtheorem{Definition}[Theorem]{Definition}

\newtheorem{Lemma}[Theorem]{Lemma}
\newtheorem{Remark}[Theorem]{Remark}
\begin{document}
\title{On the excursions of drifted Brownian motion and the successive passage times of  Brownian motion }
\author{Mario Abundo\thanks{Dipartimento di Matematica, Universit\`a  ``Tor Vergata'', via della Ricerca Scientifica, I-00133 Rome, Italy.
E-mail: \tt{abundo@mat.uniroma2.it}}
}
\date{}
\maketitle

\begin{abstract}
\noindent By using the law of the excursions of Brownian motion with drift, we find the distribution
of the $n-$th passage time of Brownian motion through a straight line $S(t)= a + bt.$
In the special case when $b = 0,$ we  extend the result to a space-time transformation of
Brownian motion.
\end{abstract}


\section{Introduction}
We consider a drifted Brownian motion of the form
\begin{equation} \label{DBM}
X(t)= x + \mu t + B_t, \ t \ge 0,
\end{equation}
where $x, \mu \in \mathbb{R}$ and $B_t$ is standard Brownian
motion (BM). When $X(t)$ is entirely positive or entirely negative on the time interval $(a,b),$
it is said that it is an excursion of drifted BM;
this means that $B_t$ remains above or
below the straight line $-x - \mu t$, for all $t \in (a,b).$
Excursions of drifted BM have interesting applications in Biology, Economics, and other applied sciences.
As an example in Economics,
if we admit that the time evolution of the gross domestic product is described by a BM with drift $\mu ,$ starting from $x$,
then
the downward and upward movement of it around its long-term growth trend (i.e. the straight line $x + \mu t),$
gives rise to an economic cycle.
These fluctuations typically involve shifts over time between periods of relatively rapid economic growth
(expansions or booms), and periods of
relative stagnation or decline (contractions or recessions) (see e.g. \cite{sch:39}).
Excursions of drifted BM are also related to the last passage time of BM through a linear boundary;
actually, last passage times of continuous martingales play an important
role in Finance, for instance, in models of  default risk (see e.g. \cite{ell:2000}, \cite{jea:2009}).
\par
When the drift $\mu$ is zero, $X(t)$ becomes BM and it is well-known that the excursions of BM have the arcsine law, namely
the probability that BM has no zeros in the time interval $(a,b)$ is given by $\frac 2 \pi \arcsin \sqrt {a/b} $ (see e.g. \cite{kleb:05}).
By using Salminen's formula for the last passage time of BM through a
linear boundary (see \cite{salm:88}), we find the law of the excursions of drifted BM, namely
the probability that
$X(t)= x + \mu t + B_t$ has no zeros in the interval $(a,b).$ From this, we
derive the distribution of the $n-$th passage time of BM
through the linear boundary $S(t)= a + bt, \ t \ge 0.$ \par\noindent
We recall
that the first-passage time of BM  through $S(t),$ when starting
from $x ,$  is defined by $ \tau _1 (x) = \inf \{ t >0: x +B_t = a + bt \} $ and
the Bachelier-Levy formula holds:
$$
P(\tau_1 (x) \le t ) = 1 - \Phi((a -x) / \sqrt t + b \sqrt t )
 + exp(-2b (a-x)) \Phi(b \sqrt t -(a-x)/ \sqrt t ),
$$
where $\Phi (y)= \int _ {- \infty } ^y \phi(z) dz ,$ with
$\phi (z) = e^ {- z^2/2}/ \sqrt {2 \pi },$
is the cumulative distribution function of the standard Gaussian variable.
If $(a -x)b > 0 ,$ then $P(\tau _1 (x) < \infty )= e^{-2b(a-x)},$ whereas, if
$(a -x)b \le 0 , \ \tau _1 (x)$  is finite with probability one and it has the following Inverse Gaussian density,
which is non-defective
(see e.g. \cite{kar:98}):
\begin{equation} \label{IGdensity}
f_ {\tau _1} (t) = f_ {\tau _1} (t |x) = \frac { d} {dt } P(\tau _1 (x) \le t) =  \frac {|a-x| } {t^ {3/2} } \ \phi \left ( \frac { a +bt -x } {\sqrt t } \right ) ,
\ t >0;
\end{equation}
moreover, if $b \neq 0,$ the expectation of $\tau _1 (x)$ is finite, being
$E(\tau _1(x)) = \frac {|a-x| } {|b| }.$ \par
The second-passage time of BM  through
$S(t),$ when starting from $x ,$ is defined by
$ \tau _2 (x) = \inf \{ t > \tau _1(x): x +B_t = a + bt \} ,$ and generally, for $n \ge 1, \  \tau _n (x) = \inf
\{ t > \tau _{n-1}(x): x +B_t = a + bt \} $ denotes the $n-$th passage
time of BM through $S(t).$
Our aim is to study its distribution. \par
The paper is organized as follows: in Section 2
we will find explicitly the distribution of the $n-$th passage
time of BM, in Section 3,
we will extend the result to space-time transformations of BM, in the special case when $b=0.$

\section{The n-th passage time of Brownian motion }
In this section we suppose that $ b \le 0 $ and $x < a ,$ or $b \ge 0 $ and $ x > a,$ so that  $P( \tau _1(x) < \infty )=1.$
First, for fixed $t>0,$ we consider the last-passage-time prior to $t$ of BM, starting from $x ,$ through the boundary
$S(t)=a + bt,$ that is:
\begin{equation} \label{LPT}
\lambda _S ^t = \sup \{ 0 \le u \le t : x + B_u = S(u) \}.
\end{equation}
The distribution of $\lambda _S ^t$ can be expressed in terms of the first-passage-time distribution of BM
through the time-reversed boundary $\widehat S(u)=S(t-u)$ (see \cite{salm:88}); in particular, we can derive from \cite{salm:88} the
following formula for the probability density, say $\psi_t (u),$  of $\lambda _S ^t :$
\begin{equation} \label{salminenEQ}
\psi _t (u)= \frac d {du} P (\lambda _S ^t \le u ) = \frac 1 { \sqrt {2 \pi u }} \exp (-b^2 u/2) \int _ {- \infty} ^{+ \infty }
\nu _ {x-a} (t-u, \widehat S )  dx, \ u \le t
\end{equation}
where:
\begin{equation} \label{nu}
\nu _x (v, \widehat S ) = \exp \{ - b(x-bt) - b^2 v /2 ) \} \frac {|x-bt|} {\sqrt { 2 \pi v^3} } \exp \{ - (x-bt)^2 /2v \}
\end{equation}
Then, the following explicit formula is obtained,  by calculation:

\begin{Lemma} \label{lemmaunico}
The probability density of $\lambda _S ^t$ is explicitly given by:
\begin{equation} \label{distrLPT}
 \psi _t (u)=  \frac { e^{ - \frac {b^2} 2 u} } { \pi \sqrt {u(t-u)}}
\left [ e^{ - \frac {b^2} 2 (t-u)} + \frac b 2 \sqrt { 2 \pi  (t-u)} \ \Big ( 2 \Phi ( b \sqrt {t-u} \ ) -1 \Big ) \right ]  , \ 0 <u <t
\end{equation}
where $\Phi (y)= \int _ {- \infty } ^y \phi(z) dz $ is the distribution function of the standard Gaussian variable. Notice that
$ \psi _t $ is independent of $a.$
\par\noindent
In particular, if $b=0,$ one gets:
\begin{equation} \label{distrLPTb=0}
\psi _t (u) = \frac 1 { \pi \sqrt {u(t-u)}} \ ,  \ 0 <u <t ,
\end{equation}
that is,  the arc-sine law with support in $(0,t).$
\end{Lemma}
\bigskip

\noindent {\it Proof.} \ By using \eqref{salminenEQ} and \eqref{nu}, we obtain:
\begin{equation} \label{equatforpsi}
 \psi _t (u)= \frac {e^{ - b^2u /2 }} { \sqrt {2 \pi u } } \cdot  \frac {e^{ - b^2 (t-u)/2 }} {t-u} \cdot J (u) ,
 \end{equation}
where
$$ J (u) = \int _ { - \infty } ^{+ \infty } \frac {e^{-b y} |y| } { { 2 \pi } \sqrt {t-u}  } e^{ -y^2/ 2 (t-u) } dy .$$
Setting $z= y / \sqrt {t-u} ,$  the integral $J$ assumes the form:
$$ J(u)= \sqrt { t-u} \int _ {- \infty } ^ {+ \infty } \frac {|z|} {\sqrt { 2 \pi } }  e^{ - (z^2 + 2z b \sqrt {t-u} )/2} dz $$
$$ = \sqrt {t-u} \ e^{ b^2 (t-u)/2 } \int _ {- \infty } ^{+ \infty } \frac {|z| } {\sqrt { 2 \pi } } \
e ^{ - (z + b \sqrt {t-u} )^2 /2} dz $$
$$ = \sqrt {t-u} \ e^{ b^2 (t-u)/2 } \left [ \int _ { - \infty } ^{b \sqrt {t-u}} \frac {b \sqrt {t-u} - w  } {\sqrt {2 \pi } } \
e^ {-w^2 /2} dw + \int _ {b \sqrt {t-u} } ^{+ \infty } \frac {w - b \sqrt {t-u}  } {\sqrt {2 \pi } } \
e^ {-w^2 /2} dw  \right ] .$$
By direct calculations, we get:
$$ J (u)= \sqrt {t-u} \ e^{ b^2 (t-u)/2 } \left [ \frac 2 {  \sqrt {2 \pi } } \ e ^ {- b^2 (t-u)/2 }
+ 2b \sqrt {t-u} \ \Phi ( b \sqrt {t-u} \ ) - b \sqrt {t-u} \ \right ].$$
Finally,  \eqref{distrLPT} soon follows, by using \eqref{equatforpsi}.

\hfill $\Box$

\begin{Remark}
For $z <t,$ the event $ \{ \lambda _S ^t \le z \} $ is nothing but the event $ \big \{ x+ B_u - S(u)$
{ \it has no zeros in the  interval} $(z,t) \big \} $.
\end{Remark}

Now, we go to consider the second-passage-time, $\tau _2 (x),$ of BM starting from $x,$ through the linear boundary $S(t)=a + bt ,$
when $ x <a $ and $b \le 0,$ the case when $b \ge 0$ and $x > a$ can be studied in a similar way.
We set $T_1(x)= \tau _1(x)$ and $T_2 (x) = \tau_2(x) - \tau _1 (x).$ We will see that
$ T_2 (x)$ is finite with probability one only if $b =0.$
Conditionally to $\tau _1 (x)=s,$ the event $ \{ \tau _2 (x) > s+t \} \ (t >0) ,$ is nothing but the event $ \{ x + B_u - S(u)$
{ \it has no zeros in the interval } $(s, s+t) \} .$ Therefore, from Lemma \ref{lemmaunico}:
$$ P( \tau _2 (x) > \tau _1 (x) +t | \tau _1 (x) =s ) =  P( \lambda _S ^{s+t} \le s ) = \int _0 ^s \psi _ {s+t} (y) dy $$
and so
\begin{equation} \label{conditionaldistrT2}
P( T_2(x) \le t | \tau _1 (x) =s ) = 1 - \int _ 0 ^s \psi _{s+t} (y) dy .
\end{equation}
Then:
$$
P( T_2 (x) \le t )= \int _ 0 ^{+ \infty } \left [ 1 - \int _ 0 ^s \psi _{s+t} (y) dy \right ] f_{\tau _1} (s) ds
$$
\begin{equation} \label{distrT2}
= 1- \int _ 0 ^{+ \infty } f_{\tau _1} (s) ds \int _ 0 ^s \psi _{s+t} (y) dy .
\end{equation}
By taking the derivative with respect to $t,$ we obtain the density of $T_2(x):$
$$
f_{T_2} (t) = - \int _ 0 ^{+ \infty } f_{\tau _1} (s) ds  \frac \partial { \partial t } \int _ 0 ^s \psi _ {s+t} (y) dy $$
\begin{equation} \label{denT2}
= \int _ 0 ^{+ \infty } f_{\tau _1} (s)  \left [ e^{ - b^2 (s+t)/2 } \frac {\sqrt s } {\pi \sqrt t (s+t) } \right ] \ ds .
\end{equation}
Notice that $ \int _ 0 ^s \psi _{s+t} (y) dy $ is decreasing in $t,$ that is,
$\frac \partial { \partial t } \int _0 ^s  \psi _ {s+t} (y) dy = \int _0 ^s \frac \partial { \partial t }  \psi _ {s+t} (y) dy <0.$
As easily seen,
$ f_{T_2} (t) \sim  const / \sqrt t,$ as $ t \rightarrow 0^+.$
Moreover, from \eqref{distrT2} it follows that, if $b \neq  0,$ the distribution of $T_2(x)$ is  defective, that is, $P( T_2(x) = + \infty ) >0.$
In fact:
$$ \lim _ { t \rightarrow + \infty } \int _ 0 ^s \psi _ {s+t} (y) dy = \frac {|b| } {\sqrt {2 \pi } }
\int _0 ^s \frac { e^{-b^2 y/2}} {\sqrt y } \ dy = 2 \ {\rm sgn} (b) \left ( \Phi ( b \sqrt s ) - \frac 1 2 \right ), $$
where ${\rm sgn} (b) =
\begin{cases}
0 & {\rm if} \  b=0 \\  |b| / b & {\rm if} \  b \neq 0
\end{cases}.$ \par\noindent
Therefore,  we get:
$$ P( T_2(x) < + \infty ) = 1- \int _ 0 ^{ + \infty } f_{\tau _1} (s) \cdot 2 {\rm sgn} (b) \left ( \Phi ( b \sqrt s ) - \frac 1 2 \right ) ds $$
$$
= 1- 2 {\rm sgn} (b) E \left ( \Phi ( b \sqrt { \tau_1(x)} \ ) - \frac 1 2 \right ),
$$
which is less than $1,$ if $b \neq 0.$ Thus, in this case:
\begin{equation} \label{PT2infty}
P( T_2(x) = + \infty ) = 2 {\rm sgn} (b) E \left ( \Phi ( b \sqrt { \tau_1(x)} \ ) - \frac 1 2 \right ) >0.
\end{equation}
Since the function $g(s)= 2 {\rm sgn} (b) \left ( \Phi ( b \sqrt { s} \ ) - \frac 1 2 \right )$ is concave, by  Jensen's inequality
we get:
$$ P(( T_2(x) = + \infty ) \le \gamma (b), $$
where:
\begin{equation} \label{gamma}
\gamma (b) =2 {\rm sgn} (b) \left ( \Phi ( b \sqrt { E( \tau _1 (x)} \ ) - \frac 1 2 \right )
= 2 {\rm sgn} (b) \left ( \Phi ( b \sqrt {  |(a-x)/b|  } \ ) - \frac 1 2 \right ) .
\end{equation}
Notice that  $ \gamma$ is an even function of $b.$ \par
On the contrary, if $b=0,$ from \eqref{PT2infty} we get that $P(T_2(x) = + \infty ) =0,$ that is, $T_2(x)$ is a proper random variable,
and
also $\tau _2 (x) = T_2(x) + \tau _1 (x)$ is finite with probability one;
precisely, by calculating the integral in \eqref{distrT2}
we have:
\begin{equation}
P( T_2(x) \le t ) = \int _0 ^{+ \infty } \frac 2 \pi \arccos \sqrt { \frac s {s+t} } \
\frac {|a-x| } {\sqrt { 2 \pi } s^{3/2} } \ e^{ - (a -x)^2 /2s } ds
\end{equation}
and
\begin{equation} \label{denT2b0}
f_{T_2} (t) = \int _ 0 ^{+ \infty } \frac {1} {\pi (s+t) \sqrt t }
\frac {|a-x| } {\sqrt { 2 \pi } s } \ e^{ - (a -x)^2 /2s } ds .
\end{equation}
In the Figure 1 we compare the plot of $P(\tau _2(x) = + \infty) = P(T _2(x) = + \infty)$
as a function of $b\le 0,$ for $a=1$ and $x=0,$ with the plot of its upper bound $\gamma (b)$ given by \eqref{gamma}.
Notice that, as it must be, $\gamma (b)$ approaches $1,$ for large negative values of $b.$

\begin{figure}
\centering
\includegraphics[height=0.33 \textheight]{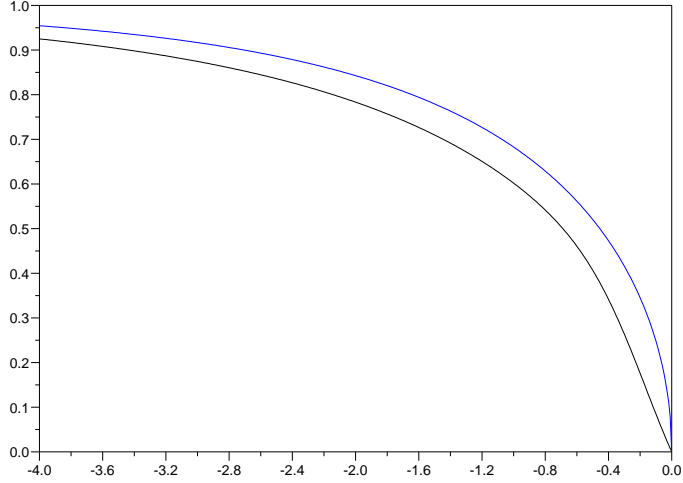}
\caption{Comparison of the shapes of $P(\tau _2(x) = + \infty)$ (lower curve)
and that of its upper bound $\gamma (b)$ (see \eqref{gamma}),  as a function of $b\le 0,$ for $a=1$ and $x=0.$
}
\end{figure}

In the Figure 2, we report the probability density of $T_2(x),$ obtained from \eqref{denT2}, by calculating numerically the integral, for
$x=0, \ a=1,$ and various values of the parameter $b \le 0.$

\begin{figure}
\centering
\includegraphics[height=0.33 \textheight]{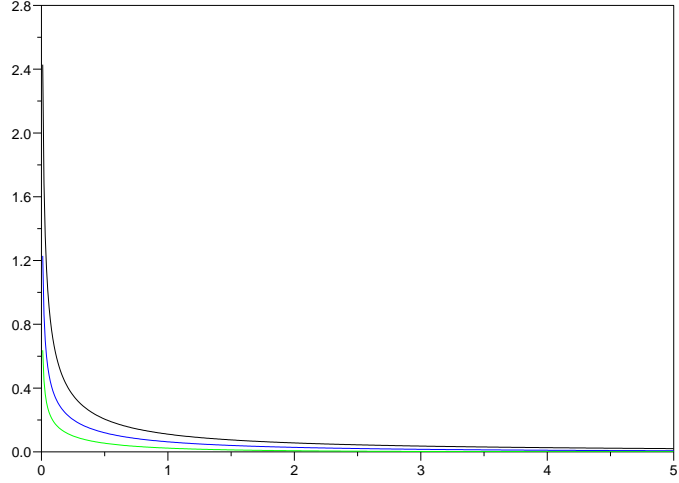}
\caption{Approximate density of $T_2(x)$ for $x=0, \ a =1$ and various values of the parameter $b;$ from top to the bottom: $b=0, \ b=-0.5, \ b=-1.$
}
\end{figure}

\indent As far as the expectation of  $\tau _2 (x)$ is concerned, it is obviously infinite for $b \neq 0 ,$ wile $E(\tau _1 (x))$ is finite.
If $b=0, \ E( \tau _1(x))$ and $E( \tau _2(x))$
are both infinite.
As for $\tau _1 (x)$ this is well-known, as for $\tau _2 (x)$ it derives from the fact that
$E(T_2(x)) = + \infty .$
Indeed, since
$$ P( T_2(x) >t )= 1- \int _0 ^{+ \infty } \frac 2 \pi \arccos \sqrt { \frac s {s+t} } \
\frac {|a-x| } {\sqrt { 2 \pi } s^{3/2} } \ e^{ - (a -x)^2 /2s } ds $$
$$= \int _0 ^{+ \infty } \frac 2 \pi \arcsin \sqrt { \frac s {s+t} } \
\frac {|a-x| } {\sqrt { 2 \pi } s^{3/2} } \ e^{ - (a -x)^2 /2s } ds ,$$
we have
$$ E(T_2 (x)) = \int _ 0 ^{+ \infty } P( T_2(x) >t ) dt =
\int _ 0 ^{ + \infty } ds  \frac 2 \pi \ \frac {|a-x| } {\sqrt { 2 \pi } s^{3/2} } \ e^{ - (a -x)^2 /2s } \int _ 0 ^{ + \infty }
\arcsin \sqrt { \frac s {s+t}} dt  \ ,$$
which is infinite,
because, as easily seen, $ \int _ 0 ^{ + \infty }
\arcsin \sqrt { \frac s {s+t}} dt = + \infty ,$ for any $s >0 .$
\par
By taking the derivative with respect to $t$ in \eqref{conditionaldistrT2},
we obtain the  density of $T_2(x)$ conditional to $\tau _1(x)=s,$ that is:
\begin{equation}
f_{T_2 | \tau _1} (t | s) = - \frac d {dt} \int _0 ^s  \psi _ {s +t } (y) dy = e^ {- b^2 (s+t)/2} \frac {\sqrt s } { \pi (s+t) \sqrt t } \ ,
\end{equation}
and, for $b=0:$
\begin{equation}
f_{T_2 | \tau _1} (t | s) = \frac {\sqrt s } { \pi (s+t) \sqrt t } .
\end{equation}
Since $ \tau _2 (x) = \tau _1 (x) + T_2 (x),$ by the convolution formula, we get the density of $\tau _2 (x):$
\begin{equation} \label{dentau2}
f_ {\tau _2} (t) = \int _ 0 ^t f_{T_2 | \tau _1} (t-s | s) f _ {\tau _1 } (s) ds  = \frac {e^ { - b^2 t/2}} { \pi t}
\int _0 ^t \frac {|a -x| } { \sqrt { 2 \pi } s \sqrt { t-s } } e^{ - (a +bs-x)^2 /2s } ds .
\end{equation}
Of course, the distribution of $\tau _2 (x)$  is defective for $b \neq 0,$ namely
$\int _ 0 ^ {+ \infty} f_ {\tau _2} (t) dt = 1 - P( \tau _2 (x)= + \infty) < 1 ,$
since $P( \tau _2 (x)= + \infty) = P( T _2 (x)= + \infty) >0.$ \par\noindent
If $b=0,$ we obtain:
\begin{equation} \label{dentau2b0}
f_ {\tau _2} (t) = \frac 1 { \pi t }
\int _0 ^t \frac {|a -x| } { \sqrt { 2 \pi } s \sqrt { t-s } } e^{ - (a -x)^2 /2s } ds ,
\end{equation}
which is non-defective. \par
In the Figure 3, we report the probability density of $\tau _2 (x)$ obtained from \eqref{dentau2} by calculating numerically the integral, for
$x=0, \ a=1,$ and various values of the parameter $b \le 0.$ Although the shapes appear to be similar to that of the inverse Gaussian density
\eqref{IGdensity}, the density
of $\tau _2 (x)$ is more concentrated around its maximum. In the Figure 4, we report the comparison between the probability density of $\tau _2 (x)$ and the inverse Gaussian density, for $a=1, b=0$ and $x=0.$

\begin{figure}
\centering
\includegraphics[height=0.33 \textheight]{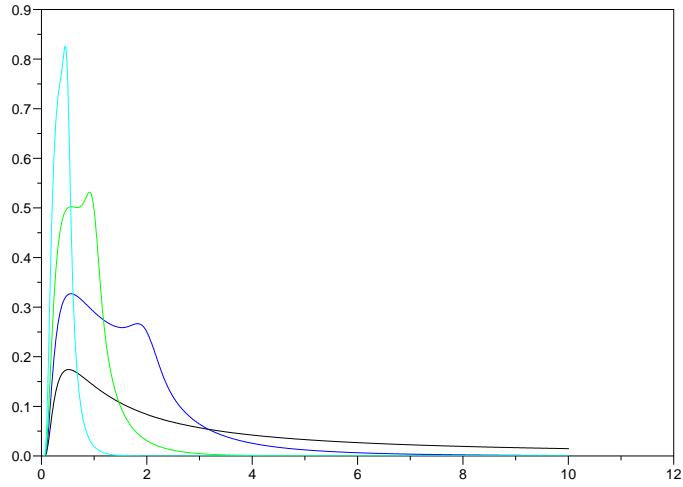}
\caption{Approximate density of $\tau _2(x)$ for $x=0, \ a =1,$ and various values of the parameter $b;$ from top to the bottom,
with respect to the peak of the curve: $b=-2, \ b=-1, \ b=-0.5, \ b=0.$
}
\end{figure}

\begin{figure}
\centering
\includegraphics[height=0.33 \textheight]{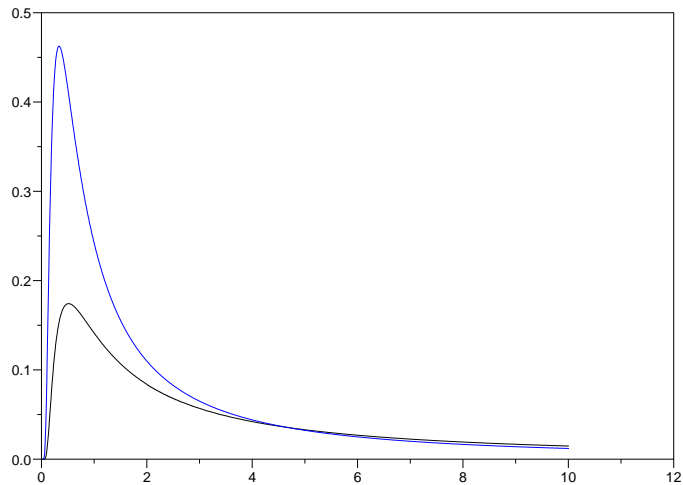}
\caption{Comparison between the probability density of $\tau _2 (x)$ (upper peak) and the inverse Gaussian density (lower peak),
for $a=1, b=0$ and $x=0.$
}
\end{figure}

\bigskip\bigskip

By reasoning in analogous manner as above,
we conclude:
\begin{Proposition}
Let be
$T_1(x) = \tau _1 (x), \ T_n (x)= \tau _n (x) - \tau _ {n-1} (x), \  n =2, \dots $  \par\noindent
Then:
\begin{equation}
P(T_1 (x) \le t ) = 2(1 - \Phi(a -x / \sqrt t  )),
\end{equation}
\begin{equation}
P( T_n (x) \le t ) = 1 - \int _ 0 ^{+ \infty } f_{\tau _ {n-1} } (s) ds \int _ 0 ^s \psi _{s+t} (y) dy, \ n= 2, \dots
\end{equation}
Moreover, the density of $\tau _n (x)$ is:
\begin{equation}
f_ {\tau _n} (t) = \int _ 0 ^t f_{T_n | \tau_{n-1}} (t-s | s) f _ {\tau _{n-1} } (s) ds ,
\end{equation}
where $f_{\tau _ {n-1} }$ and $f_ {T_n | \tau _ {n-1} }$ can be calculated inductively, in a similar way, as done for $f_{\tau _ {2} }$
and $f_{T_2 | \tau _1 } .$ \par\noindent
If $b=0, \ T_1(x), \ T_2 (x), ....$ are finite with probability one.
\end{Proposition}

\hfill $\Box$

\begin{Remark}
The expression for $P(T_1(x) \le t)$ is nothing but the  Bachelier-Levy formula, written for $b=0.$
\end{Remark}

\section{The n-th passage time of space-time transformations of Brownian motion}
The techniques of the previous section can be applied to get also results for time-changed BM. In fact, let
be
\begin{equation} \label{timechangedBM}
Z(t)= z + B(\rho (t)),
\end{equation}
where $\rho (t) \ge 0$ is an increasing, differentiable function of $t \ge 0 ,$ with $\rho (0)=0.$
Such kind of diffusion process $Z$ is a special case of Gauss-Markov process (see \cite{abundo:smj13}, \cite{din:aap01}, \cite{ric:smj08}, \cite{nob:smj06}); in particular the
form \eqref{timechangedBM} is taken by certain integrated Gauss-Markov processes (see \cite{abundo:smj15}).
Now, for $n \ge 1$ denote again by $\tau _n (z)$ the successive passage times of $Z(t)$ through the constant barrier $S=a,$ and
by $T_n (z)= \tau _n(z) - \tau _{n-1} (z)$ the inter-passage times; then, for $z <a,$ we have
$ \tau _1 (z)= \inf \{ t >0: z + B(\rho (t)) =a \}  $ and so $\rho (\tau _1 (z)) = \inf \{ s>0: z + B_s =a \} := \tau ^B _1 (z),$ where
the superscript $B$ refers to BM. Therefore, $\tau_1(z) = \rho ^ {-1} (\tau ^B _1 (z)),$ and in analogous way, for $n \ge 2 , $ we get
$\rho (\tau _n (z)) = \inf \{ s > \rho (\tau  _ {n-1} (z)) : z + B_s) =a \} :=
 \tau ^B _ {n} (z).$ In conclusion, we have:
$$ \tau _n (z)= \rho ^ {-1} ( \tau ^B _n (z)) \ {\rm and } \  T_n(z) = \rho ^ {-1} ( \tau ^B _n(z))- \rho ^ {-1} ( \tau ^B _{n-1}(z)), \ n \ge 1 ,$$
where $\tau ^B _n (z)$ is the n-th passage time of BM, starting from z, through $a;$
thus, the calculations of the distributions of $\tau _n (z)$ and $ T_n (z)$ are reduced to those of $\tau ^B _n (z)$ and
$ T^B_n (z).$ \par\noindent
An analogous study concerning the successive spike (i.e. passage) times $\tau _n $ of a Gauss-Markov  process $Z$
through a constant threshold $S,$
was developed  in \cite{pir:2015}, for a non homogeneous Leaky Integrate-and-Fire (LIF) neuronal model, in which
the membrane potential of the neuron, $Z,$ is instantaneously reset to its initial value, every time $S$ is attained.
\par
The approach considered in the present paper can be also applied to
one-dimensional diffusions which can be reduced to BM by a space transformation; let
$Z(t)$ be the the solution of the stochastic differential equation (SDE):

\begin{equation} \label{eqdiffu}
d Z(t) = \mu (Z(t)) dt + \sigma (Z(t)) d B_t \ , \  Z(0) = z ,
\end{equation}
where the coefficients $\mu (z)$ and $\sigma (z)$ are regular enough functions (see e.g. \cite{abundo:saa13}),
so that a unique strong solution exists.
We consider the following: \bigskip

\begin{Definition}
We say that $Z$
is conjugated to BM if there exists
an increasing  function $v$ with $v(0)=0,$ such that
 $Z(t) = v^ {-1} (B_t + v(z)),$ for any $t \ge 0.$
\end{Definition}
\bigskip
Examples of diffusions conjugated to BM are the following (for more, see e.g. \cite{abundo:saa13}): \par\noindent
{\bf (i)}  \ (Feller process or Cox-Ingersoll-Ross (CIR) model) \par\noindent
the solution of the SDE
$$ d Z(t) = \frac 1 4  dt + \sqrt {Z(t) \vee 0} \ dB_t \ , Z(0)=z , $$
which is conjugated to BM via the function $v(z) =2 \sqrt z ,$ i.e. $Z(t) = \frac 1 4 (B_t +2 \sqrt { z} )^2;$  \par\noindent
{\bf (ii)} (Wright \& Fisher-like process)\par\noindent
the solution of the SDE
$$d Z(t) = \left ( \frac 1 4 - \frac 1 2  Z(t) \right ) dt + \sqrt {Z(t)(1- Z(t)) \vee 0} \ dB_t \ , Z(0) = z \in [0,1] ,$$
which is conjugated to BM
via the function
$v(x) = 2 \arcsin \sqrt {z},$ i.e. $Z(t) = \sin ^2 ( B_t /2 + \arcsin \sqrt {z}).$ 
\bigskip

Let $Z$ be conjugated to BM, via the function $v,$ and denote, as always, by $\tau  _{n} (z) $ the n-th passage time of $Z$ through the constant
barrier $S=a,$  and by $T _{n} (z)$ the inter-passage time, with the condition that $Z(0)=z;$ as easily seen, since now $\rho (t)=t,$ we get:
$$ \tau   _{n} (z) =  \tau ^{B,a'} _{n} (z') \ {\rm and } \  T _ n(z) =  \tau ^{B,a'} _{n} (z')-  \tau ^{B,a'} _{n-1} (z'), \ n \ge 1 ,$$
where  $a'=v(a), \ z'=v(z)$ and $\tau ^{B,a'} _{n} (z')$ denotes the $n-$th passage time of BM through the barrier $a',$ when starting from $z'.$ Again,
the calculations of the distributions of $\tau  _n (z)$ and $ T_n (z)$ are reduced to those concerning BM. \par
This also works, for certain moving boundaries $S,$ e.g. for Geometric Brownian motion, that is
the solution of the SDE:
$$d Z(t) = r  Z(t) dt + \sigma Z(t) d B_t , \ Z(0)= z >0  , $$
where $ r$ and $\sigma$ are positive constant. This is a well-known equation in the framework of Mathematical Finance, since
it describes the time evolution of a stock price $Z;$  its explicit solution is
$Z(t)= z e ^ { \mu t } e ^ { \sigma B_t} ,$
where $\mu = r - \sigma ^2 /2 ,$ so $\ln Z / \sigma $ turns out to be BM with drift $\mu ,$ starting from $\ln z / \sigma .$
If we consider the moving barrier $S(t)= e ^ {\sigma S_0 + \mu ' t },$ then the successive passages, $\tau _n (z),$ of $Z$ through $S$ are
reduced to the successive passages, $ \tau _n ^B (z'),$ of BM through  the linear  boundary 
$S_0 + ( \mu ' - \mu )t/ \sigma ,$  with the condition that the starting point is $z'=\frac {\ln z } \sigma .$ In fact:
$$ \tau   _{n} (z) =  \rho ^ {-1} ( \tau _{n} ^B (z')) \ {\rm and } \  T _ n(z) =  \rho ^ {-1} ( \tau  _{n} ^B (z'))-
 \rho ^ {-1} ( \tau  _{n-1} ^B (z')).$$
The same considerations apply to the Ornstein-Uhlenbeck process, which is solution
of the SDE:
$$ d Z(t) = - \mu  Z (t) dt + \sigma dB_t, \  Z(0)= z ,$$
where $\mu, \sigma$ are positive constants. By using a time--change (see e.g. \cite{abundo:stapro12}),
the explicit solution assumes the form
$  Z(t)= e^{- \mu t } \left ( z + B( \rho(t)) \right ),$  where
$ \rho (t) =  \frac { \sigma ^2 } { 2 \mu } \left ( e ^ {2 \mu t} -1 \right).$
Let us consider the moving barrier $S(t)=S_0 e^{ - \mu t},$ with $S_0 > z;$ the first-passage time of $ Z (t)$ through $S(t)$ is
$  \tau _{S(t)}   =  \inf \{ t>0: z+ B (\rho (t)) = S_0 \} ,$ and so
$ \rho (  \tau_{S(t)}) =  \inf \{ u >0: z + B_u =  S_0 \} .$
Therefore, similarly to
the case of Geometric BM, the successive passages of $Z$ through $S$ are
reduced to those of BM  through  the constant  boundary $S_0.$

\end{document}